\input amstex
\documentstyle{amsppt}

\document

\topmatter
\title
Horospherical Cauchy-Radon transform on compact symmetric spaces
\endtitle

\author Simon Gindikin \endauthor

\address Departm. of Math., Hill Center, Rutgers University,
110 Frelinghysen Road, Piscataway, NJ 08854-8019
\endaddress
\email gindikin\@math.rutgers.edu \endemail \abstract Harmonic
analysis on noncompact Riemannian symmetric spaces is in a sense
equivalent to the theory of the horospherical transform. There are
no horospheres on compact symmetric spaces, but we define a
complex version of horospherical transform which plays a similar
role for the harmonic analysis on them.\endabstract
\endtopmatter
In the horospherical transform we integrate test functions along
horospheres or, in other words, we integrate these functions with
$\delta$-functions with supports on horospheres. The idea of the
complex horospherical transform \cite {Gi00} is to replace the
$\delta$-function by the Cauchy kernel: more precisely, we take
Cauchy kernels with singularities along complex horospheres
without real points. Such a transform can work when the real
transform either does not exist or has too big a kernel \cite
{Gi00, Gi01, Gi04}. We define such a transform here on compact
symmetric spaces. Earlier we saw \cite {Gi04} that it already
gives a new interesting integral transform for the sphere $S^n$.

Let $X=G/K$ be a compact Riemannian symmetric space; $G$ here is a
simply connected compact semisimple Lie group; $K$ is its
connected involutive subgroup. Let $\frak g=\frak k +\frak n$ be
the corresponding involutive decomposition in the Lie algebras.
Let $X_n=G_n/K$ be the dual noncompact Riemannian symmetric space
and $\frak g_n=\frak k+i\frak n$ be the corresponding
decomposition in Lie algebras. Let $\frak a$ be a Cartan subspace
in $i\frak n$, $\Sigma_+$ be a system of positive (restricted)
roots; $(\beta_1,\dots \beta_l)$ be a basis of unmultipliable
roots and $(\mu_1,\dots, \mu_l)$ be the dual system of fundamental
weights ($\frac{\langle \mu_i,\beta_j\rangle}{\langle
\beta_j,\beta_j\rangle}=\delta_{ij}$). Then
$$\Lambda =\{\mu=\lambda _1 \mu_1+\cdots+ \lambda_l; \lambda_j\in
\Bbb Z^+\}$$ is the set of highest (restricted) weights of
spherical representations corresponding to the symmetric space
$X$.

Let $A= \exp(\frak a)$, $A_I=\exp(i\frak a)$ be the corresponding
commutative subgroups in $G_n, G$ correspondingly. Let $G_n=KAN$
be the Iwasawa decomposition, where $N$ is the maximal unipotent
subgroup in $G$. Let $G_\Bbb C^0=K_\Bbb C A_\Bbb C N_\Bbb C$ be
the holomorphic complexification of this decomposition in $G_\Bbb
C$. For a group $H$ we denote through $H_\Bbb C$ its
complexification. The set $G_\Bbb C^0$ is a Zariski open set in
$G_\Bbb C$. Let us consider the projection $a(g)$ of $G_\Bbb C^0$
on the factor $A_\Bbb C$. We can consider this function also as a
function on the symmetric Stein space $X_\Bbb C=G_\Bbb C/K_\Bbb
C$.

For $a\in A_\Bbb C$ we define

$$a_j=\exp(\mu_j(\log a)), \quad a^\mu=\prod _{1\leq j\leq l} (a_j)^{\lambda_j}, \mu\in
\Lambda.$$ For noncompact symmetric spaces $X_n$ zonal spherical
functions can be obtained by the averaging of $a(x)^\mu$ along
$K$-orbits (Harish-Chandra's integral representation). Clerc \cite
{Cl88} suggested in the case of compact symmetric spaces $X$,
using the possibility of the holomorphic extension of zonal
spherical functions on $X_\Bbb C$, to write the analogue of
Harish-Chandra representation on $X_\Bbb C$. This representation
can be restricted on the compact group $G$ and the compact space
$X$ but there is involved the complex function $a^\mu (x)$. For
the zonal spherical function $\phi_\mu(x)$ with the center at the
point $x_0\in X$ (corresponding to a fixed isotropy subgroup $K$)
we have
$$\phi_\mu(x;x_0)=\int _K a^\mu(xk) \nu(dk),x\in X_\Bbb C$$ where we integrate on the invariant
measure on $K$. Here and everywhere below we integrate on a
compact homogeneous manifold $Y$ on the invariant measure
$\nu(dy)$ normalized by the condition $\int_Y \nu(dy)=1$. The
integrand is well-defined for almost all $k\in K$. Clerc used this
representations for the computation of the asymptotic of zonal
spherical functions on compact symmetric spaces similar to the
famous Harish-Chandra's computation for noncompact spaces. The
construction of the complex horospherical transform on $X$ in this
note continues this application of complex geometry to the
harmonic analysis on compact symmetric spaces.

The functions $a^\mu(x), x\in X_\Bbb C,$ are elements of eigen
subspaces of the invariant differential operators which are
$N_\Bbb C$-invariant. They sometimes are called zonal
horospherical functions. Let us extend them for different choices
of $A,N$. Let
$$h_\mu(x, g)=a^\mu(xg), x\in X_\Bbb C, g\in G_\Bbb C.$$ This function
will not change under right multiplications on elements of $N_\Bbb
C$ as well as on elements of $M_\Bbb C$ (as usual, $M$ is the
centralizer of $A$ in $K$). Let $\Xi=G_\Bbb C/M_\Bbb C N_\Bbb C$.
We will realize its points as classes $\zeta=(gM_\Bbb C N_\Bbb
C)$. Then the function $h_\mu$ will depend only on the class
$\zeta(g)\in \Xi.$ Correspondingly, we will denote it as
$h_\mu(x,\zeta),\quad x\in X_\Bbb C, \zeta \in \Xi.$

 We will call $\Xi$ the horospherical
manifold. This name is connected with the possibility to realize
points of $\Xi$ as the horospheres on $X_\Bbb C$: nondegenerate
orbits of the subgroup $N_\Bbb C$ and its $G_\Bbb C$-translations.

Let $F=G_\Bbb C/P_\Bbb C=G_\Bbb C/M_\Bbb C A_\Bbb C N_\Bbb C\cong
G/MA$ be the corresponding flag manifold. There is a natural
fibering $\Xi \rightarrow F$ with fibers $A_\Bbb C$. The group
$G_\Bbb C$ acts on $\Xi$ by the "left multiplications"
($\zeta\mapsto g\cdot \zeta $) and $A_\Bbb C$  by the "right
multiplications" ($\zeta \mapsto \zeta\cdot a))$ and these actions
commutate. The existence of the right action of $A_\Bbb C$ is a
crucial circumstance for the horospherical transform. We have
$$h_\mu(x, \zeta\cdot a)=a^\mu h_\mu(x,\zeta), x \in X_\Bbb C,
\zeta \in \Xi, a\in A_\Bbb C.$$ Choosing an initial point
$\zeta_0$ we can parameterize points $\zeta \in\Xi$ as pairs
$\zeta=(a,u), a\in A_\Bbb C, u\in F$ such that the right
multiplications on elements of $A_\Bbb C$ act on the
$a$-components of $\zeta=(a,u)$. We can define the horospheres
$\Omega(\zeta), \zeta\in \Xi,$ on $X_\Bbb C$ by the equations
$$h_j(\cdot, \zeta)=1, \quad 1\leq j\leq l,$$ where
$h_j=h_{\mu_j}$. Of course, then $h_\mu(x,\zeta)=1$ for all $\mu$.

We have on $X$ \cite {Cl88,Corollary, p.427}: $$|a_j(x)|\geq 1,
1\geq j\leq l$$ and the conditions $|a_j(x)|=1$ distinguish the
image of the subgroup $A_I$ in $X$ ($a_j(x)=1$ for the initial
point $x_0$). Let $A_+$ be the semigroup $\{|a_j|< 1, 1\leq j\leq
l\}$ in $A_\Bbb C$ and we define the domain $\Xi_+$ in $\Xi$ as
$\{\zeta=(a,u); a\in A_+\}$ (the $(G\cdot \zeta_0\cdot A_+ )$-
orbit of the initial point $\zeta_0$). So the domain $\Xi_+$ is
fibering over the flag domain $F$ with the fibers $A_+$. We can
reformulate the definition of $\Xi_+$: $$ \Xi_+=\{\zeta \in
\Xi;h_j(x,\zeta)<1\, \text {for all}\, x\in X \,\text
{and}\,j\}.$$

In other words, the domain $\Xi_+$ parameterizes horospheres
$\Omega(\zeta)$ on $X_\Bbb C$ without real points (they do not
intersect $X$). The points of the boundary $\zeta\in\partial
\Xi_+$ correspond to the horospheres which intersect $X$ in one
point. It means that the boundary admits the fibering over the
symmetric space $X$ with the fiber $\Xi(x)$ over $x\in X$ defined
by the equations
$$h_j(x,\zeta)=1,1\leq j\leq l,$$ (or $a(x \cdot \zeta) =e$).
The fiber $\Xi_x$ parameterizes the horospheres which intersect
$X$ in the point $x$. The fibers are isomorphic to $K/M$.

Now we are ready to give the basic definition. For $f\in C^\infty
(X)$ we define its horospherical Cauchy-Radon transform as
$$\hat f(\zeta)=\int _X \frac {f(x)} {\prod _{1\leq j\leq l}
(1-h_j(x,\zeta))} \nu(dx),\zeta \in \Xi_+.$$ For $\zeta\in \Xi_+$
the integrand has no singularities on $X$ and is holomorphic on
$\zeta$. Hence the function $\hat f$ is holomorphic in the domain
$\Xi_+$. The singular set of the integrand on $X_\Bbb C$ has a
very simple structure and has as the edge exactly the horosphere
$\Omega(\zeta)$. For $\zeta \in
\partial \Xi_+$ the integral exists in the distribution sense and
gives the regular boundary values of $\hat f$.

As usual in integral geometry, the principal problem is the
problem of the inversion of this horospherical transform. Let
$D_j$ be the partial logarithmic derivative in the direction of
the basic root $\beta_j$ on $A_\Bbb C$ and its transfer on $\Xi$
relative to the right multiplications on $A_\Bbb C$. Then $D_j
h^\mu(x,\zeta)=\mu_j h^\mu(x,\zeta).$ Let for any root $\alpha\in
\Sigma_+$:
$$D(\alpha)=1+\frac {\langle \alpha, D\rangle} {\langle\alpha, \rho\rangle}$$
where $\rho$ is the half-sum of the positive roots and $$ \Cal
L=\prod_{\alpha\in \Sigma_+} D(\alpha).$$ \proclaim {Theorem}
There is an inversion formula for the horospherical Cauchy-Radon
transform

$$f(x)=\int _{\Xi(x)} \Cal L\hat
f(\zeta) \nu(d\zeta).$$ \endproclaim Here we integrate the
boundary values of the holomorphic function $\hat f$ over the
normalized invariant measure on $K/M$ .

This result is a simple corollary of the Plancherel formula on
compact symmetric spaces \cite {Sh77,H94}. We define the spherical
Fourier  transform in "complex" form: $$\tilde f(\zeta;
\mu)=\int_X f(x)h^\mu (x,\zeta) \nu(dx),\quad \mu \in \Lambda,
\zeta\in \Xi.$$ The function $\tilde f(\zeta; \mu)$ is holomorphic
on $\zeta\in \Xi$ and satisfies the condition of homogeneity:
$$\hat f(\zeta \cdot a; \mu) =a^\mu(\hat f(\zeta; \mu).$$ In other
words, $\tilde f(\zeta;\mu)$ is the section of the line bundle on
$F$ corresponding to $\mu$. So we have maps in the spaces of
irreducible spherical representations in the flag realization.

Similarly, to noncompact symmetric spaces, the spherical Fourier
transform and horospherical Cauchy-Radon transform are connected
by the commutative Fourier transform, this time, by the discrete
one. Since the right multiplications on elements $a\in A_\Bbb
C$($\zeta \mapsto \zeta \cdot a$) commutate with the action of $G$
we can decompose $\hat f(\zeta)$ in the Fourier series relative to
this action of Abelian compact group $A_I$. The components (which
are invariant relative to the action of $G$) are just $\tilde
f(\zeta;\mu)$:

$$\hat f(\zeta)=\sum_{\mu\in \Lambda} \tilde f(\zeta ;\mu).$$
It is the result of the direct decomposition the kernel in the
definition of the horospherical transform. The inversion formula
for the spherical Fourier transforms equivalent to the Plancherel
formula on $X$ is
$$f(x)=\sum_{\mu\in \Lambda} d(\mu)\int _{\Xi(x)} \tilde
f(\zeta;\mu)\nu(d\zeta),$$ where for the dimension we have Weyl's
formula
$$d(\mu)=\prod_{\alpha\in \Sigma_+} \frac {\langle \mu
+\rho,\alpha\rangle}{\langle \rho,\alpha\rangle}.$$ It follows
from the inversion formula on $X$ \cite {Sh77,H94} since $$\int
_{\Xi(x)} \tilde f(\zeta,\mu)\nu(d(\zeta)=\int_X
f(y)\phi(y;x)\nu(dy).$$ Here $\phi(\cdot,x)$ is the zonal
spherical function with the center $x$ and it follows from Clerk's
representation for zonal spherical functions. Now the inversion
formula is the consequence of the fact that the multiplication of
the Fourier coefficients on the polynomial $d(\mu)$ corresponds to
the action of the differential operator $\Cal L$.

In this note we obtain the results on horospherical transform as a
reformulation of known facts of the harmonic analysis on compact
symmetric spaces. Our aim here is to emphasize that formulas which
we receive for the complex horospherical transform on compact
symmetric spaces are completely similar to the formulas for the
horospherical transform on noncompact Riemannian symmetric spaces.
In our next publication we will obtain them by direct methods of
integral geometry and it gives a new approach to harmonic analysis
on compact homogeneous manifolds.

Let us consider a simplest example when $X$ is the sphere
$S^n\subset \Bbb R^{n+1}$ \cite {Gi04}:$$\Delta(x)=x\cdot
x=(x_1)^2+ \cdots +(x_{n+1})^2=1.$$ Let $\Bbb C S^n$ be its
complexification defined in $\Bbb C^{n+1}$ by the equation $\Delta
(z)=1$. Then $\Xi$ can be realize as the cone
$$\Delta(\zeta)=0, \zeta \in \Bbb C^{n+1}, \zeta \ne0,$$ and the
flag manifold $F$ as its projectivization in $P\Bbb C^n$. The
horospheres $\Omega(\zeta),\zeta \in \Xi,$ can be defined as the
intersection of $S^n$ by the isotropic hyperplanes
$$\zeta\cdot z=\zeta_1 z_1+\cdots +\zeta_{n+1} z_{n+1}=1.$$
The parameters $\zeta$ of horospheres without real points give the
domain $\Xi_+$:
$$\Delta (\xi)=\Delta(\eta)<1, \zeta=\xi+i\eta.$$
The fibers $\Xi(x)$ in the fibering of $\partial \Xi_+$ over $S^n$
are $\{\zeta=x+i\eta, x \cdot \eta=0,\eta\cdot \eta=1\}$. Finally
the differential operator in the inversion formula is $\Cal
L=(1+cD)^{n-1}$, where $D$ is the logarithmic differentiation
along generators in the cone $\Xi$ and $c$ is a normalized
constant.

\remark  {Remarks} 1.The domain $\Xi$ covers the compact complex
manifold $F$. On $F$ all holomorphic functions are constant, but
on $\Xi_+$ there is a lot of holomorphic functions (they separate
points of $\Xi_+$. The domain $\Xi_+$ is not holomorphically
complete since $\Xi$ does not have this property. The extension
reduces to the extension of the fibers $A_\Bbb C$. If $F=\Bbb
CP^n$ and $z_0,\dots, z_n$ are homogeneous coordinates, then
$\Xi_+$ is the complex ball without the point: $1>|z_0|^2+\cdots
+|z_n|^2>0$ and all holomorphic functions in $\Xi_+$ extend on the
ball. In the example of the sphere, which we considered earlier,
the natural extension of $\Xi_+$ has the singularity.

It is simple to construct a Hardy type norm in $H(\Xi_+)$
corresponding to the action of $A_+$ and then we will have the
model of spherical series of representations for $X$.

2. The horospherical transform and its inversion connect functions
on real manifold $X$ and holomorphic functions in $\Xi_+$. The
natural understanding of this situation lies in the consideration
of hyperfunctions on $X$ as $\bar \partial$-cohomology classes of
$H^{(m)}(X_\Bbb C \backslash X)$, where $m=n-l, n=\dim X, l=
\operatorname{rank} X$ (the dimension of horospheres). Then we can
interpret the horospherical transform as a version of the Penrose
transform: we integrate cohomology along complex horospheres in
$X_\Bbb C \backslash X$. A complication is that horospheres are
not cycles and we need to be careful when integrating cohomology.
It would be interesting to define natural subspaces of cohomology
in $X_\Bbb C \backslash X$ and holomorphic functions in $\Xi_+$
which can be intertwined. Intermediate Paley-Wiener theorems can
be also interesting.

To receive the cohomological interpretation of the inverse
horospherical transform we need to appeal to holomorphic language
for $\bar \partial$-cohomology from \cite {EGW95}. For each $u\in
F$ we take the union $D(u)$ of "parallel" horospheres
parameterizing by points in $\Xi_+$ lying over $u$. We have the
covering of $X_\Bbb C \backslash X$ by the Stein domains $D(u),
u\in F$. This covering satisfies the conditions from \cite {EGW95}
and we can construct Dolbeault cohomology using the complex of
holomorphic forms $\varphi (z|u,du)$ on $F$ depending
holomorphically on parameters $z\in D(u)$. We can extend the
integrand in the inversion formula as such kind of a form.

There is another possibility interesting from point of view of
complex analysis: to consider the horospherical Cauchy-Radon
transform as the intertwining operator between holomorphic
functions on whole manifolds $X_\Bbb C$ and $\Xi$. We will discuss
this operator in another paper.

3.When we want to construct the inverse horospherical transform we
seek an integral operator on the space of all in a sense functions
on the manifold of horospheres. In our case they are holomorphic
functions. The restrictions of this operator on subspaces of
sections of line bundles on $F$, corresponding to irreducible
representations, give eigen functions of invariant differential
functions (spherical polynomials). It is a Poisson's type
integral. For the sphere such formulas go back to Maxwell. It
would be interesting to deliberate the analogy with noncompact
case. In a sense the flag manifold $F$ plays the role of the
"complex" boundary of $X$.

\endremark

\enddocument